\documentstyle[amsfonts,12pt]{article}
\newcommand{\f}{{\cal O}}
\newcommand{\lon}{\longrightarrow}
\newcommand{\rar}{\rightarrow}
\newcommand{\hook}{\hookrightarrow}
\newcommand{\pr}{pr}
\newcommand{\OM}{\Omega^1 M}
\newcommand{\CP}{{\Bbb C} {\Bbb P}}
\newcommand{\rank}{\mbox{rank}}

\newtheorem{theorem}{Theorem}
\newtheorem{proposition}[theorem]{Proposition}
\newtheorem{lemma}[theorem]{Lemma}
\newtheorem{definition}[theorem]{Definition}

\begin{document}

\title{Affine connections on involutive $G$-structures}

\author{Sergey A.\ Merkulov\\
\small    Department of Pure Mathematics, Glasgow University\\
\small  15 University Gardens, Glasgow G12 8QW, UK}

\date{ }
\maketitle
\sloppy

\paragraph{0.\ Introduction.}
 An affine connection  is one of the central objects in differential
geometry.  One of its most informative characteristics  is the
(restricted) holonomy group which is defined, up to a conjugation,
as a subgroup of $GL(T_t M)$ consisting of all
automorphisms of the tangent space $T_t M$ at a point $t\in M$
induced by parallel translations along the $t$-based contractible
loops in $M$. The list of groups which can be  holonomies of affine
connections is dissappointingly dull ---  according to Hano and
Ozeki \cite{HO}, {\em any}\, closed subgroup of a general linear
group can be realized in this way. The situation, however, is very
different in the subclass of affine connections with {\em zero
torsion}. Long ago, Berger \cite{Berger} presented a very restricted
list of possible
irreducibly acting holonomies of torsion-free  affine connections.
His list was complete in the part of {\em metric}\, connections (and
later much
work has been done to refine this "metric" part of his list, see, e.g.,
\cite{Bryant2} and references cited therein), while
the situation with holonomies of {\em non-metric}\,  torsion-free affine
connections was and remains very unclear. One of the
results that will be discussed  in this paper asserts that  {\em any}\,
torsion-free holomorphic affine
connection with irreducibly acting holonomy group can, in principle,
be constructed by twistor methods. Another result reveals a new
 natural subclass of affine connections with   {\em very little
torsion}\, which shares  with the class of torsion-free affine
connections two basic properties --- the list of irreducibly acting
holonomy groups of affine
connections in this subclass is  very restricted and the links with the
twistor theory are again very strong.

The purpose of this paper is to explain the key elements of the above
mentioned twistor constructions without indulging in rather lengthy
proofs.
We work throughout in the category of complex manifolds,
holomorphic affine connections, etc., though many results
can be easily adapted to the
real analytic case along the lines explained in \cite{Manin}.

\paragraph{1.\ Irreducible $G$-structures.}
When studying an affine connection $\nabla$  with the irreducibly acting
holonomy group $G$, it is  suitable  to work with the
associated $G$-structure. In this section we recall some notions
of the theory of $G$-structures.

Let $M$ be an $m$-dimensional complex manifold and ${{\cal L}}^*M$ the
holomorphic coframe bundle $\pi: {{\cal L}}^*M \rar M$ whose fibres
$\pi^{-1}(t)$ consist of all ${\Bbb C}$-linear isomorphisms
$e: {\Bbb C}^{m} \rar \Omega^1_t M$, where $\Omega^1_t M$ is the
cotangent space at $t\in M$. The space ${{\cal L}}^*M$ is a principle
right $GL(m,{\Bbb C})$-bundle with the right action given by
$R_g(e)= e\circ g$.
If $G$ is a closed subgroup of $GL(m,{\Bbb C})$, then a (holomorphic)
$G$-structure on $M$ is a principal subbundle ${\cal G}$ of
${{\cal L}}^*M$ with the group $G$. It is clear that there is a
one-to-one correspondence between the set of $G$-structures on $M$
and the set of holomorphic sections $\sigma$
of the quotient bundle $\tilde{\pi}: {{\cal L}}^*M/G \rar M$ whose
typical fibre is isomorphic to $GL(m,{\Bbb C})/ G$.
A $G$-structure on $M$ is called {\em locally
flat}\, if there exits a coordinate patch in the neighbourhood of
each point $t\in M$ such that in the associated canonical
trivialization of ${{\cal L}}^*M/G$ over this patch the section
$\sigma$ is represented
by a constant $GL(m,{\Bbb C})/G$-valued function.
A G-structure is called {\em k-flat}\, if,
for each $t\in M$, the $k$-jet of the associated section $\sigma$ of
${{\cal L}}^*M/G$ at $t$ is isomorphic to the $k$-jet of some
locally flat section of ${{\cal L}}^*M/G$. It is not difficult to
show that a $G$-structure admits a torsion-free affine connection if
and only if it is 1-flat (cf.\
\cite{Bryant2}). A $G$-structure on $M$ is called {\em irreducible}\,
if the action of $G$ on ${\Bbb C}^{m}$ leaves no non-zero invariant
subspaces.

Given an affine connection $\nabla$ on a connected
simply connected complex manifold $M$ with the irreducibly acting
holonomy group $G$, the associated irreducible  $G$-structure
${\cal G}_{\nabla}\subset{{\cal L}}^*M$ can be constructed as follows.
Define two points $u$ and $v$ of ${{\cal L}}^*M$ to be equivalent,
$u\sim v$, if there is a
holomorphic path $\gamma$ in $M$ from $\pi(u)$ to $\pi(v)$ such that
$u=P_{\gamma}(v)$, where $P_{\gamma}:\Omega^1_{\pi(v)}M\rar
\Omega^1_{\pi(u)}M$ is the parallel transport along $\gamma$.
Then ${\cal G}_{\nabla}$ can be defined, up to an isomorphism,
as $\left\{u\in {{\cal L}}^*M \mid u\sim v\right\}$ for some coframe
$v$. The $G$-structure ${{\cal G}}_{\nabla}$ is
the smallest subbundle of ${{\cal L}}^*M$ which is invariant under
$\nabla$-parallel translations.

It will be shown later that for any holomorphic irreducible
$G$-structure ${\cal G}\rar M$ there is associated an analytic family
of compact isotropic submanifolds $\{X_t \hook Y\mid t\in M\}$ of a
certain complex contact manifold $Y$ which encodes much information
about ${\cal G}$. To explain this correspondence in more detail, we
first digress in the next two sections to the Kodaira \cite{Kodaira}
deformation theory of compact
complex submanifolds and to its particular generalization
studied in \cite{Me1}.

\paragraph{2. Kodaira relative deformation theory.}
Let $Y$ and $M$ be complex manifolds and let $\pi_{1}: Y\times
M\lon Y$ and $\pi_{2} : Y\times M \lon M$ be  natural projections.
An analytic family of compact submanifolds of the complex manifold
$Y$ with the parameter space $M$ is a complex submanifold
$F\hookrightarrow Y\times M$ such that the
restriction of the projection $\pi_{2}$ on $F$ is a proper regular map
(regularity means that the rank of the differential of $\nu :=
\left.\pi_{2}\right|_{F}: F\lon M$
is equal to $\dim M$ at every point). The parameter space $M$ is
called a {\em Kodaira moduli space}. Thus the family $F$ has
a double fibration structure
$$
Y \stackrel{\mu}{\longleftarrow} F \stackrel{\nu}{\lon} M
$$
where $\mu := \left.\pi_{1}\right|_{F}$.  For each $t\in M$ we say
that the compact complex submanifold $X_{t} = \mu\circ\nu^{-1}(t)
\hookrightarrow Y$
belongs to the family $F$. Sometimes we use a more explicit
notation $\{X_t\hook Y\mid t\in M\}$ to denote an analytic
family $F$ of compact submanifolds.

If $F\hook Y\times M$ is an analytic family of compact submanifolds,
then, for any $t\in M$, there is a natural linear map \cite{Kodaira},
$$
k_t : T_{t}M \lon H^{0}(\nu^{-1}(t), N_{\nu^{-1}(t)\mid F})
\stackrel{\mu_{*}}{\lon}
H^{0}(X_{t}, N_{X_{t}\mid Y}),
$$
which is a composition of the natural lift of a tangent vector at $t$
to a global section of the normal bundle of the submanifold $\nu^{-1}(t)
\hookrightarrow F$ with the Jacobian of $\mu$ (here the symbol $N_{A|B}$
stands for the normal bundle of a submanifold $A\hookrightarrow B$). An
analytic family $F\hookrightarrow Y\times M$ of compact submanifolds is
called {\em complete} if the  map $k_{t}$ is an isomorphism for all
$t\in M$
which in particular implies that $\dim M = h^0(X_{t}, N_{X_{t}\mid Y})$.

In 1962 Kodaira \cite{Kodaira} proved the following existence
theorem: {\em if $X \hookrightarrow Y$ is a compact complex submanifold
with normal bundle $N$ such that $H^{1}(X,N) = 0$, then $X$ belongs to
a complete analytic family $F\hook Y\times M$ of compact submanifolds
of $Y$.}

\paragraph{3. Deformations of compact complex Legendre submanifolds
of complex contact manifolds.}
In this section we shall be interested in the following specialisation
(which will eventually turn out to be a generalisation)
of the  Kodaira relative deformation problem: the initial data is a pair
$X\hook Y$ consisting of a compact complex Legendre submanifold $X$ of
a complex contact manifold $Y$ and the object of study is the set, $M$,
of all holomorphic deformations of $X$ inside $Y$ which remain Legendre.
First, we recall  some standard notions,  then give a better formulation
of the problem, and finally present its solution.

Let $Y$ be a complex $(2n + 1)$-dimensional manifold. A complex contact
structure on $Y$ is a rank $2n$ holomorphic subbundle $D\subset TY$ of
the holomorphic tangent bundle to $Y$ such that the Frobenius form
\begin{eqnarray*}
\Phi: D \times D & \longrightarrow & TY/D\\
(v,w) & \longrightarrow & [v,w]\bmod D
\end{eqnarray*}
is non-degenerate. Define the contact line bundle $L$ by the exact
sequence
$$
0\lon D \lon TY  \stackrel{\theta}{\lon} L \lon 0.
$$
One can easily verify that maximal non-degeneracy of the distribution
$D$ is equivalent to the fact that the above defined "twisted" 1-form
$\theta \in H^0(Y, L\otimes \Omega^1 M)$ satisfies the condition
$\theta\wedge (d\theta)^n \neq 0$.
A complex submanifold $X\hook Y$ is called {\em isotropic}\, if
$TX\subset D$.
An isotropic submanifold of maximal possible dimension $n$ is called
{\em Legendre}. In this paper we shall be primarily interested in
compact Legendre submanifolds.
The normal bundle $N_{X\mid Y}$ of any Legendre submanifold
$X\hookrightarrow Y$ is
isomorphic to $J^1 L_X$ \cite{L2}, where $L_X = \left.L\right|_X$.
Therefore, $N_{X\mid Y}$ fits into the exact sequence
$$
0\lon \Omega^1 X\otimes L_{X} \lon N_{X\mid Y} \stackrel{pr}{\lon}
L_{X} \lon 0 . \label{ext}
$$

Let $Y$ be a complex contact manifold. An analytic family
$F\hook Y\times M$ of compact submanifolds of $Y$ is called
an {\em analytic family of compact Legendre submanifolds} if, for any
point $t\in M$,
the corresponding subset $X_{t}:= \mu\circ\nu^{-1}(t) \hookrightarrow Y$
is a Legendre submanifold.
The parameter space $M$ is called a {\em Legendre moduli space}.

Let $F\hook Y\times M$ be an analytic family of compact Legendre
submanifolds. According to Kodaira~\cite{Kodaira}, there is a natural
linear map $k_t : T_{t}M \lon H^{0}(X_{t}, N_{X_{t}\mid Y})$.
We say that the family $F$ is {\em complete at a point}\, $t\in M$ if
the composition
$$
s_{t}:  T_{t}M \stackrel{k_{t}}{\lon} H^{0}(X_{t}, N_{X_t\mid Y})
\stackrel{\pr}{\lon} H^{0}(X_{t}, L_{X_{t}})
$$
provides an isomorphism between the tangent space to $M$ at the point
$t$ and the vector space of global sections of the contact line bundle
over $X_{t}$. One of the motivations behind this definition is the
fact \cite{Me1} that  an analytic family of compact Legendre
submanifolds \{$X_{t}\hook Y \mid t \in M$\}  which is complete at a
point $t_{0}\in M$ is also  {\em maximal}\, at
 $t_{0}$ in the sense that, for any other analytic family
of compact Legendre submanifolds  \{$X_{\tilde{t}}\hook Y \mid \tilde{t}
\in \tilde{M}$\}
such that $X_{t_{0}} = X_{\tilde{t}_0}$
for a point $\tilde{t}_{0}\in \tilde{M}$, there exists a neighbourhood
$\tilde{U}\subset \tilde{M}$ of $\tilde{t}_0$ and a holomorphic map
$f: \tilde{U}\lon M$ such that $f(\tilde{t}_{0}) = t_{0}$ and
$X_{f(\tilde{t}')} = X_{\tilde{t}'}$ for each $\tilde{t}'\in \tilde{U}$.
An analytic family  $F\hook Y\times M$
is called {\em complete}\, if it is complete at each point of the
Legendre moduli space $M$. In this case $M$ is also called complete.

The following result \cite{Me1} reveals a simple condition for the
existence of complete Legendre moduli spaces.

\begin{theorem}\label{leg}
Let $X$ be a compact complex Legendre submanifold of a complex contact
manifold
$(Y,L)$. If $H^{1}(X,L_{X}) = 0$, then there exists a complete
analytic family  of compact Legendre submanifolds $F\hook Y\times M$
containing $X$. This family is maximal and $\dim M = h^{0}(X,L_{X})$.
\end{theorem}

Let $X$ be a complex manifold and $L_X$ a line bundle on $X$. There is a
natural "evaluation" map
$H^0(X,L_X)\otimes \f_X \lon J^1 L_X$
whose dualization gives rise to the canonical map
$$
L_X\otimes S^{k+1} (J^1 L_X)^* \lon L_X\otimes S^{k} (J^1 L_X)^*
\otimes \left[H^0(X,L_X)\right]^*
$$
which in turn gives rise to the map of cohomology groups
$$
H^1\left(X,L_X\otimes S^{k+1} (J^1 L_X)^*\right) \stackrel{\phi}{\lon}
H^1\left(X,L_X\otimes S^{k} (J^1 L_X)^*\right)
\otimes \left[H^0(X,L_X)\right]^* .
$$
For future reference, we define a vector subspace
$$
\tilde{H}^1\left(X,L_X\otimes S^{k+1} (J^1 L_X)^*\right) := \ker \phi
\, \subset \, H^1\left(X,L_X\otimes S^{k+1} (J^1 L_X)^*\right).
$$

\paragraph{4. $G$-structures induced on Legendre moduli spaces
of generalized flag varieties.}
Recall that a generalised flag variety $X$ is a compact simply
connected homogeneous K\"{a}hler manifold \cite{BE}. Any such a manifold
is of the form $X=G/P$, where $G$ is a complex semisimple Lie
group and $P\subset G$ a fixed parabolic subgroup.
Assume that such
an $X$ is embedded as a Legendre submanifold into a complex contact
manifold $(Y,L)$ with contact line bundle $L$ such that
$L_X := \left.L\right|_X$ is very ample. Then the
Bott-Borel-Weil theorem and the fact that any holomorphic line bundle on
$X$ is homogeneous imply that $H^1(X,L_X)=0$.
Therefore, by Theorem~\ref{leg},  there exists a complete  analytic
family of compact Legendre
 submanifolds $\{X_t\hook Y\mid t\in M\}$, i.e.\ the initial data
"$X\hook Y$" give rise to  a new complex manifold $M$ which,
as the following result shows, comes equipped with a rich
geometric structure.

\begin{theorem}{\em \cite{Me1}}\label{Main}
Let $X$ be a generalised flag variety embedded as a Legendre submanifold
into a complex contact manifold $Y$ with contact line bundle $L$ such that
$L_X$ is very ample on $X$. Then
\begin{description}
\item[(i)] There exists a complete analytic family $F\hook Y\times M$ of
compact  Legendre submanifolds with moduli space
$M$ being an $h^0(X,L_X)$-dimensional complex manifold. For each
$t\in M$, the associated Legendre submanifold $X_t$ is isomorphic to $X$.

\item[(ii)] The Legendre moduli space $M$ comes equipped with an induced
irreducible $G$-structure, ${{\cal G}}_{ind} \rar M$,
with $G$ isomorphic to the connected component of the identity of
the group of all global biholomorphisms $\phi: L_X
\rightarrow L_X$ which commute with the projection $\pi: L_X\rightarrow X$.
The Lie algebra  of $G$  is isomorphic to
$H^0\left(X, L_X\otimes (J^1 L_X)^*\right)$.

\item[(iii)] If ${{\cal G}}_{ind}$ is $k$-flat, $k\geq 0$, then the
obstruction for
${{\cal G}}_{ind}$ to be $(k+1)$-flat is given by a tensor field on $M$
whose value at each $t\in M$ is represented by a cohomology class
$\rho_t^{[k+1]}\in \tilde{H}^1\left(X_t, L_{X_t}
\otimes S^{k+2}(J^1 L_{X_t})^*\right)$.

\item[(iv)] If ${{\cal G}}_{ind}$ is 1-flat, then the bundle of all
torsion-free connections in ${{\cal G}}_{ind}$ has as the typical
fiber an affine space
modeled on $H^0\left(X, L_{X}\otimes S^2(J^1 L_{X})^*\right)$.
\end{description}
\end{theorem}
\vspace{3 mm}

\noindent{\sc Remark}. Theorem~\ref{Main} is actually valid for
a larger class of compact complex manifolds $X$ than the class of
generalized flag varieties --- the only vital assumptions are \cite{Me1}
that $X$ is rigid and the cohomology groups $H^1(X,\f_X)$
and $H^1(X,L_X)$ vanish.
\vspace{3 mm}

The geometric meaning of cohomology classes
$\rho_t^{[k+1]}\in \tilde{H}^1\left(X_t, L_{X_t}
\otimes S^{k+2}(J^1 L_{X_t})^*\right)$ of Theorem~\ref{Main}(iii)
is very simple --- they compare to $(k+2)$th order the germ of the
Legendre
embedding $X_t\hook Y$ with the "flat" model, $X_t\hook J^1 L_{X_t}$,
where the ambient contact manifold is just the total space of the vector
bundle $J^1 L_{X_t}$ together with its canonical contact structure
and the Legendre submanifold $X_t$ is realised as a zero section of
$J^1 L_{X_t} \rar {X_t}$. Therefore, the cohomology class
$\rho_t^{[k]}$ can be called the $k$th Legendre jet of $X_t$ in $Y$.
 Then it is natural to call a Legendre submanifold
$X_t\hook Y$\, $k$-{\em flat}\, if $\rho_t^{[k]}=0$.
With this terminology, the item (iii) of Theorem~\ref{Main} acquires
a rather symmetric form: {\em the induced $G$-structure on the moduli
space $M$ of a complete analytic family of compact Legendre submanifolds
 is $k$-flat if and only if the family consists of $k$-flat Legendre
submanifolds}.

This general construction can be illustrated by three well
known examples which were among the motivations behind the present work
(in fact the list of examples can be made much larger ---
Theorem~\ref{Main} has been checked for all "classical" torsion-free
geometries as well as for a large class of locally symmetric structures).
The first example  is a "generic" $GL(m,{\Bbb C})$-structure on an
$m$-dimensional manifold $M$. The associated twistorial
data $X\hook Y$ is easy to describe: the complex contact manifold $Y$ is
the  projectivized cotangent bundle ${\Bbb P}(\OM)$ with its natural
contact
structure while $X =\CP^{m-1}$ is just a fiber of the projection
${\Bbb P}(\OM)\rar M$. The corresponding complete family
$\{X_t\hook Y\mid t\in M\}$ is the set of all fibres of this fibration.
Since $L_X=\f(1)$ and $J^1 L_X={\Bbb C}^m\otimes\f_X$,
we have $H^1\left(X,L_X\otimes S^{k+2}((J^1L_X)^*)\right)=0$
for all $k\geq 0$ which confirms the well-known
fact that any $GL(m,{\Bbb C})$-structures on an
$m$-dimensional manifold are locally flat.

The second example \cite{L1}
is a pair $X\hook Y$ consisting of an $n$-quadric $Q_n$
embedded into a $(2n+1)$-dimensional contact manifold $(Y,L)$ with
$\left. L\right|_X\simeq i^*\f_{{{\Bbb C} {\Bbb P}}^{n+1}}(1)$,
$i: Q_n\rightarrow {{\Bbb C} {\Bbb P}}^{n+1}$ being a standard projective
realisation of
$Q_n$. It is easy to check that in this case \mbox{$H^0(X,L_X\otimes
(J^1 L_X)^{*})$} is precisely the conformal algebra implying that the
associated $(n+2)$-dimensional Legendre moduli space $M$ comes equipped
canonically with a conformal structure.
Since $H^1\left(X, L_X\otimes S^2(J^1 L_X)^*\right)=0$, the induced
conformal structure must be torsion-free in agreement with the
classical result of differential geometry. Easy calculations show
that the vector space $H^1\left(X, L_X\otimes S^3(J^1 L_X)^*\right)$
is exactly the subspace of $TM\otimes\OM \otimes \Omega^2M$ consisting
of tensors with Weyl curvature symmetries. Thus Theorem~\ref{Main}(iii)
implies the well-known
Schouten conformal flatness criterion. Since $H^0\left(X, L_X\otimes S^2
(J^1 L_X)^*\right)$ is isomorphic to the typical fibre of $\OM$, the set
of all torsion-free affine connections preserving the induced conformal
structure is the affine space modeled on $H^0(M,\OM)$, again in agreement
with the classical result.

The third example
is Bryant's \cite{Br} relative deformation problem $X\hookrightarrow Y$
with $X$ being a rational Legendre curve ${{\Bbb C}{\Bbb P}}^1$ in a
complex contact 3-fold $(Y,L)$ with $L_X=\f (3)$. Calculating
$H^0(X,L_X\otimes (J^1 L_X)^{*})$, one easily concludes that the induced
$G$-structure on the associated
4-dimensional Legendre moduli space is exactly an   exotic
$G_3$-structure
which has been studied by Bryant in his search for irreducibly acting
holonomy groups of torsion-free affine connections which are missing
in the Berger list \cite{Berger} (the missing holonomies are called
{\em exotic}). Since
$H^1\left(X, L_X\otimes S^2(J^1 L_X)^*\right)=0$, Theorem~\ref{Main}(iii)
says
that the induced $G_3$-structure  is torsion-free in accordance
with \cite{Br}. Since $H^0\left(X, L_X\otimes S^2(J^1 L_X)^*\right)=0$,
${\cal G}_{ind}$ admits a unique torsion-free affine connection $\nabla$.
The cohomology class $\rho_t^{[2]}\in H^1\left(X, L_X\otimes
S^3(J^1 L_X)^*\right)$
from Theorem~\ref{Main}(iii) is exactly the curvature tensor of $\nabla$.

How large is the family of $G$-structures which can be constructed by
twistor methods of Theorem~\ref{Main}? As the following result
\cite{Me1} shows, in the category of irreducible 1-flat $G$-structures
this class as large as one could wish.

\begin{theorem}\label{Main2}
\begin{description}
\item[(i)] Let $H$ be one of the following representations:
(a) $Spin(2n+1,{\Bbb C})$ acting on ${\Bbb C}^{2^n}$,   $n\geq 3$; (b)
$Sp(2n,{\Bbb C})$ acting on ${\Bbb C}^{2n}$,  $n\geq 2$; (c)
$G_2$ acting on ${\Bbb C}^7$. Suppose that
$G\subset GL(m,{\Bbb C})$
is a connected semisimple Lie subgroup whose decomposition into a
locally direct product of simple groups contains $H$. If ${\cal G}$
is any irreducible 1-flat $G\cdot {\Bbb C}^*$-structure on  an
$m$-dimensional manifold $M$, then there exists a complex contact
manifold $(Y,L)$ and a generalised flag variety $X$ embedded into $Y$
as a Legendre submanifold with $L_X$ being very ample, such that, at
least locally, $M$ is canonically isomorphic to the associated Legendre
moduli space and ${\cal G} \subset {\cal G}_{ind}$. In
particular, when $G=H$ one has in the case
(a) $X= SO(2n+2,{\Bbb C})/U(n+1)$ and ${\cal G}_{ind}$ is a
$Spin(2n+2,{\Bbb C})\cdot {\Bbb C}^*$-structure; in the case (b)
$X={\Bbb C} {\Bbb P}^{2n-1}$ and ${\cal G}_{ind}$ is a
$GL(2n,{\Bbb C})$-structure;
and in the case (c) $X=Q_5$ and ${\cal G}_{ind}$ is a
$CO(7,{\Bbb C})$-structure.

\item[(ii)]
Let $G\subset GL(m,{\Bbb C})$ be an arbitrary connected
semisimple Lie
subgroup whose decomposition into a locally direct product of
simple groups does not contain any of the groups $H$ considered in (i).
If ${\cal G}$ is any irreducible 1-flat $G\cdot {\Bbb C}^*$-structure
on an $m$-dimensional manifold $M$, then there exists a complex
contact manifold
$(Y,L)$ and a Legendre submanifold $X\hook Y$ with $X=G/P$ for some
parabolic
subgroup $P\subset G$ and with $L_X$ being very ample, such that, at least
locally, $M$ is canonically isomorphic to the associated Legendre moduli
space and ${\cal G} = {\cal G}_{ind}$.
\end{description}
\end{theorem}

The conclusion is that there are very few irreducible $G$-structures
which can
{\em not}\, be constructed by twistor methods discussed in this paper.
It is also worth pointing out that Theorem~\ref{Main}(iii) gives rise
to a new and rather effective machinery to search for exotic
holonomies.
The new results in this direction will be discussed elsewhere  ---
here we only note  that the claimed efficiency of the twistor
technique is
largerly due to the simple observation  that the key cohomology
groups  ${H}^1\left(X_t, L_{X_t}\otimes S^{2}(J^1 L_{X_t})^*\right)$
and ${H}^1\left(X_t, L_{X_t}\otimes S^{3}(J^1 L_{X_t})^*\right)$,
which provide
us with the full information about torsion and curvature tensors,
can be computed by a {\em combination}\, of the representation theory
methods (such as Bott-Borel-Weil theorem) and the methods
of complex analysis. In some  important cases it is  even enough to use
the complex analysis methods only.


\paragraph{5. Torsion-free affine connections.}
Let $F\hook Y\times M$ be a complete analytic family of compact
Legendre submanifolds. Any point
$t$ in $M$ is represented thus by a compact complex Legendre
submanifold $X_t$. The first floors  of the two towers of
infinitesimal neighbourhoods of the analytic spaces $t\hook M$ and
$X_t\hook Y$  are related to each other via the
isomorphism $T_tM = H^0(X_t, L_{X_t})$. What happens at the
second floors of these two towers? If $J_t\subset \f_M$ is the ideal
of holomorphic functions which vanish
at $t\in M$, then the tangent space $T_tM$   is isomorphic to
$\left(J_t/J_t^2\right)^*$. Define a second order tangent bundle,
$T_t^{[2]}M$, at the point $t$ as $\left(J_t/J_t^3\right)^*$. Then,
evidently,
$T_t^{[2]}M$ fits  into an exact sequence of complex vector spaces
\begin{equation}
0\lon T_tM \lon T_t^{[2]}M \lon S^2(T_tM) \lon 0 \label{mu}
\end{equation}
For each $t\in M$ there
exists a holomorphic vector bundle, $\Delta_{X_t}^{[2]}$, on the
associated Legendre submanifold $X_t\hook Y$ such that there are an exact
sequence of locally free sheaves
\begin{equation}
0\lon L_{X_t} \stackrel{\alpha}{\lon} \Delta_{X_t}^{[2]} \lon S^2(J^1
L_{X_t}) \lon 0
\label{mumu}
\end{equation}
and a commutative diagram
\begin{equation}
\begin{array}{rccccccccl}
0 & \lon & T_t M & \lon & T^{[2]}_t M & \lon & S^2(T_t M) & \lon & 0 \\
&& \Big\downarrow && \Big\downarrow && \Big\downarrow && \\
0 & \lon & H^0(X_t,L_{X_t}) & \lon &
H^0\left(X_t,\Delta_{X_t}^{[2]}\right) &
\lon & H^0\left(X_t, S^2(J^1 L_{X_t})\right) & \lon & 0
\end{array} \label{mumumu}
\end{equation}
which extends the canonical isomorphism $T_tM \rightarrow
H^0(X_t, L_{X_t})$ to
{\em second}\, order infinitesimal neighbourhoods of
$t\hookrightarrow M$
and $X_t\hookrightarrow Y$. For the details of the construction of
$\Delta_{X_t}^{[2]}$ we refer the interested reader to \cite{Me1}.
In this paper we  need only  to know  that this bundle exists and
has the stated
properties. The extension (\ref{mumu}) defines a cohomology class
$$
\rho_t^{[1]}\in\mbox{Ext}^1_{\f_{X_t}}
\left(S^2(J^1 L_{X_t}),L_{X_t}\right)
=H^1(X_t, L_{X_t}\otimes S^2(J^1 L_{X_t})^*).
$$
 This is exactly the
class of Theorem~\ref{Main}(iii) which is the obstruction to
1-flatness of $X_t$ in $Y$. Therefore, if $X_t$ is 1-flat, then
extension (\ref{mumu}) splits, i.e.\  there exists a morphism
$\beta: \Delta_{X_t}^{[2]}
\rightarrow L_{X_t}$ such that $\beta\circ\alpha = id$.
Any such a morphism induces via the commutative diagram (\ref{mumumu})
an associated splitting of the exact sequence
(\ref{mu}) which is  equivalent to a torsion-free affine connection at
$t\in M$.  A torsion-free connection
on the Legendre moduli space $M$ which arises at each $t\in M$ from a
splitting of the
extension~(\ref{mumu}) is called an {\em induced}\, connection.
Now we can formulate the main theorem about torsion-free affine
connections.

\begin{theorem}{\em \cite{Me1}}\label{conn}
Let $\nabla$ be a holomorphic torsion-free affine connection on a
complex manifold $M$ with irreducibly acting reductive holonomy
group $G$. Then there exists a complex contact manifold $(Y,L)$ and
a 1-flat Legendre submanifold  $X\hook Y$ with $X=G_s/P$ for some
parabolic subgroup $P$ of the semisimple
factor $G_s$ of $G$ and with $L_X$ being very ample, such that, at
least locally, $M$ is canonically isomorphic to the associated
Legendre moduli space
and $\nabla$ is an induced torsion-free affine connection in
${\cal G}_{ind}$.
\end{theorem}

The conclusion is that any holomorphic torsion-free affine connection
with irreducibly acting holonomy group can, in principle, be
constructed by twistor methods.


\paragraph{6. From Kodaira to Legendre moduli spaces and back.}
In this subsection we first show that {\em any}\, complete Kodaira moduli
space can be interpreted as a complete Legendre moduli space and then use
this fact to prove a proposition about canonically induced geometric
structures on Kodaira moduli spaces.

If $X\hookrightarrow Y$ is a complex submanifold, there is an exact
sequence of vector bundles
$$
0 \lon N_{X\mid Y}^{*} \lon \left.\Omega^{1}Y\right|_{X} \lon
\Omega^{1}X \lon 0,
$$
which induces a natural embedding, ${\Bbb P}(N_{X\mid Y}^{*})
\hook{\Bbb P}(\Omega^{1}Y)$, of total spaces of the associated
projectivised bundles. The
manifold $\hat{Y} = {\Bbb P} (\Omega^{1}Y)$ carries a natural contact
structure
such that the constructed embedding $\hat{X} =
{\Bbb P}(N_{X\mid Y}^{*})\hook
\hat{Y}$ is a Legendre one \cite{Ar}.
Indeed, the contact distribution $D \subset T\hat{Y}$ at each point
$\hat{y}\in\hat{Y}$ consists of those tangent vectors
$V_{\hat{y}}\in T_{\hat{y}}\hat{Y}$
which satisfy the equation $<\hat{y}, \tau_{*}(V_{\hat{y}})> = 0$, where
$\tau:
\hat{Y} \lon Y$ is a natural projection and the angular brackets denote
the pairing of 1-forms and vectors at $\tau(\hat{y})\in Y$.
Since the submanifold
$\hat{X}\subset \hat{Y}$ consists precisely of those projective classes
of 1-forms in $\left.\Omega^{1}Y\right|_{X}$ which vanish when
restricted on
$TX$, we conclude that $T\hat{X}\subset \left.D\right|_{\hat{X}}$.
One may check that this association
$$
\begin{array}{ccc}
\mbox{Kodaira moduli space} & \lon & \mbox{Legendre moduli space}\\
\{X_t\hook Y\mid t\in M\}   & \lon &
\left\{\hat{X}_{t} := {\Bbb P}( N_{X_t\mid Y}^{*})
\hookrightarrow \hat{Y}:={\Bbb P}(\Omega^{1}Y)\ \mid t\in M\right\}
\end{array}
$$
preserves completeness while changing its meaning, i.e.\ a {\em complete}\,
Kodaira family of compact complex submanifolds is mapped into a
{\em complete}\, family of compact complex Legendre submanifolds (which is
usually not complete in the Kodaira sense).

The contact line bundle $L$ on $\hat{Y}$ is just the dual of the
tautological
line bundle $\f_{\hat{Y}}(-1)$.  Simplifying the notations,
$N:= N_{X\mid Y}$
and $\hat{N}:= N_{\hat{X}\mid \hat{Y}}$, we write down the following
commutative diagram which explains how $\hat{N}$
is related to $\rho^*(N)$ and $L$
$$
\begin{array}{rccccccccl}
\vspace{2 mm} & & 0 && 0 &&&& \\
\vspace{2 mm}&& \downarrow && \downarrow && && \\
\vspace{2 mm}&& \rho^*(\Omega^1X)\otimes L_{\hat{X}} & = &
\rho^*(\Omega^1X)\otimes L_{\hat{X}} &&&& \\
\vspace{2 mm}&& \downarrow && \downarrow && && \\
\vspace{2 mm}0 & \lon & \Omega^1\hat{X}\otimes L_{\hat{X}} &
\lon & \hat{N} & \lon & L_{\hat{X}} & \lon & 0 \\
\vspace{2 mm}&& \downarrow && \downarrow && || && \\
\vspace{2 mm}0 & \lon & \Omega^1_{\rho}\otimes L_{\hat{X}} &
\lon & \rho^*(N) & \lon & L_{\hat{X}} & \lon & 0 \\
\vspace{2 mm}&& \downarrow && \downarrow && && \\
&& 0 && 0 &&
\end{array}
$$
Here $L_{\hat{X}} = \left. L \right|_{\hat{X}}$, $\rho$ is a natural
projection $\hat{X} \rar X$, and $\Omega^1_{\rho}$
is the bundle of $\rho$-vertical 1-forms, i.e.\ the dual of
$T_{\rho} = ker: T\hat{X}\rar TX$. Using this diagram it is not
hard to show that there is a long exact sequence of cohomology groups
$$
\begin{array}{rccccccl}
0&\rar& H^0\left(X, N\otimes S^2(N^*)\right)& \rar&
H^0\left(\hat{X},L_{\hat{X}}\otimes S^2(\hat{N}^*)\right)&\rar&
H^0(X, N^*\otimes TX)& \rar \\
&\rar & H^1\left(X, N\otimes S^2(N^*)\right)& \rar&
H^1\left(\hat{X},L_{\hat{X}}\otimes S^2(\hat{N}^*)\right)&\rar&
H^1(X, N^*\otimes TX)&\rar \ldots
\end{array}
$$

\begin{proposition}\label{ko}
Let $X\hook Y$ be a compact complex rigid submanifold with rigid normal
bundle $N$ such that $H^1(X,N)=0$ and let $M$ be the associated Kodaira
moduli space. If
\begin{equation}
H^1\left(X,N\otimes S^2(N^*)\right)
= H^1\left(X,N^*\otimes TX\right) = 0, \label{po}
\end{equation}
 then the associated Kodaira
moduli space $M$ comes equipped with an induced 1-flat $G$-structure
with the Lie algebra $g$ of $G$ being characterized
by the following exact sequence of Lie algebras
$$
0\lon H^0(X,N\otimes N^*) \lon g  \lon H^0(X,TX)\lon 0
$$
\end{proposition}
{\em Proof}. By Kodaira theorem \cite{Kodaira}, there is a complete
family, $\{X_t\hook Y\mid t\in M\}$, of compact complex submanifolds and
hence the associated complete family,
$\{\hat{X}_t\hook\hat{Y}\mid t\in M\}$,
of compact complex Legendre submanifolds with $\hat{X}_t =
{\Bbb P}(N_t^*)$ and $\hat{Y}= {\Bbb P}(\OM)$.
Equations~(\ref{po}) together with the above long exact sequence of
cohomology
groups imply the vanishing of
$H^1\left(\hat{X},L_{\hat{X}}\otimes S^2(\hat{N}^*)\right)$.
This together with Theorem~\ref{Main}(iii) and the subsequent Remark
implies in turn that the induced $G$-structure on $M$ is 1-flat.
The final statement about the Lie algebra of $G$  follows from
Theorem~\ref{Main}(ii),
the exact sequence
$$
0 \lon L_{\hat{X}}\otimes\rho^*(N^*) \lon L_{\hat{X}}\otimes \hat{N}^*
\lon \rho^*(TX) \lon 0
$$
and the fact that $H^0(\hat{X},L_{\hat{X}}\otimes\rho^*(N^*)) =
H^0(X,N\otimes N^*)$. $\Box$
\vspace{3 mm}

If, for example, $X$ is the projective line $\CP^1$ embedded into a complex
manifold $Y$ with normal bundle $N={\Bbb C}^{2k}\otimes \f(1)$, $k\geq 1$,
then by Proposition~\ref{ko} the associated Kodaira moduli space $M$, which
is a $4k$-dimensional complex manifold, comes equipped canonically with
a complexified quaternionic structure, in accordance with \cite{P,PP}.
For other results on geometric structures induced on
Kodaira moduli spaces we refer to \cite{Me4,MP}.


\paragraph{7. Involutive $G$-structures.}
Let $M$ be a complex $m$-dimensional manifold and let
${\cal G}\subset {\cal L}^*M$ be an irreducible holomorphic
$G$-structure with reductive $G$ (hence $G$ must be isomorphic to
$G_s$ or $G_s\times {\Bbb C}^*$ for some
semisimple Lie group  $G_s\subset GL(m,{\Bbb C})$). Since ${\cal G}$ is
irreducible, there is a naturally associated subbundle $\tilde{F}\subset
\OM$ whose typical fiber is the cone in ${\Bbb C}^m$ defined as
the $G_s$-orbit of the line spanned by a highest weight vector.
Denote $\tilde{\cal F} = \tilde{F}\setminus 0_{\tilde{F}}$, where
$0_{\tilde{F}}$ is the "zero" section of
$\tilde{p}: \tilde{F} \rar M$ whose
value at each  $t\in M$ is the vertex of the cone $\tilde{p}^{-1}(t)$.
The quotient bundle $\nu: {\cal F} = \tilde{\cal F}/{{\Bbb C}^*} \lon M$
is then a subbundle of the projectivized cotangent bundle ${\Bbb P}(\OM)$
whose fibres $X_t$ are isomorphic to the generalised flag variety
$G_s/P$, where $P$ is the parabolic subgroup of $G_s$ which preserves the
highest weight vector in ${\Bbb C}^m$ up to a scale.
The total space of the cotangent bundle $\OM$ has a canonical holomorphic
symplectic 2-form $\omega$ which makes the sheaf of holomorphic
functions on
$\OM$ into a sheaf of Lie algebras via the Poisson bracket
$\{f,g\}= \omega^{-1}(df,dg)$.

\begin{definition}
An irreducible $G$-structure ${\cal G} \rar M$ is called
involutive if $\tilde{\cal F}$ is a  coisotropic submanifold of
the symplectic manifold $\OM\setminus 0_{\OM}$.
\end{definition}

The first motivation behind this definition is the following

\begin{lemma}
Every irreducible 1-flat $G$-structure is involutive.
\end{lemma}
{\em Proof}. It is well known that a submanifold of a symplectic
manifold is isotropic if and only if the associated ideal sheaf
is the sheaf of Lie algebras relative to the Poisson bracket.
This condition obviously holds for a locally flat $G$-structure.
Since the
Poisson bracket involves only a 1st order differential operator,
this condition must also be satisfied for a 1-flat $G$-structure.
$\Box$
\vspace{3 mm}

The pullback, $i^*{\omega}$, of the symplectic form $\omega$ from
$\OM\setminus 0_{\OM}$ to its submanifold $i:
\tilde{\cal F}\lon \OM \setminus 0_{\OM}$ defines a distribution
${\cal D}\subset T\tilde{\cal F}$ as the kernel of the natural
"lowering of
indices" map $T\tilde{\cal F} \stackrel{\lrcorner\,\omega}{\lon}
\Omega^1\tilde{\cal F}$, i.e.\
${\cal D}_e = \left\{V\in T_e\tilde{\cal F}: V\lrcorner\,
i^*\omega = 0\right\}$
at each point $e\in \tilde{\cal F}$. Using the fact that
$d(i^*\omega) = i^*d\omega = 0$, one can show that this distribution is
integrable and thus
defines  a foliation of $\tilde{\cal F}$ by holomorphic leaves.
We shall assume from now on that the space of leaves, $\tilde{Y}$, is a
complex manifold.
This assumption imposes no restrictions on the local structure of $M$.
The fact that the Lie derivative, ${\cal L}_V i^*\omega =
V\lrcorner\, i*d\omega + d(V\lrcorner\, i^*\omega) = 0$, vanishes for any
vector field
$V$ tangent to the leaves implies that $i^*\omega$ is the pullback
relative to the canonical projection $\tilde{\mu}: \tilde{\cal F}
\rar \tilde{Y}$ of a closed 2-form $\tilde{\omega}$ on $\tilde{Y}$. It is
easy  to check  that $\tilde{\omega}$ is non-degenerate which means
that $(\tilde{Y}, \tilde{\omega})$ is a symplectic manifold.
There is a natural action of ${\Bbb C}^*$ on $\tilde{\cal F}$
which leaves ${\cal D}$ invariant and thus induces an action of
${\Bbb C}^*$ on
$\tilde{Y}$. The quotient $Y := \tilde{Y}/{\Bbb C}^*$ is an odd
dimensional complex manifold which has a double fibration structure
$$
Y \stackrel{\mu}{\longleftarrow} {\cal F}=\tilde{\cal F}/{\Bbb C}^*
\stackrel{\nu}{\lon} M
$$
and thus contains an analytic family of  compact submanifolds
$\left\{X_t =\mu\circ\nu^{-1}(t)\hook Y\mid t\in M\right\}$ with
$X_t = \tilde{X}_t/{\Bbb C}^*\simeq G_s/P$. Next, inverting a
well-known procedure of
symplectivisation of a contact manifold \cite{Ar}, it is not hard to show
that $Y$ has a complex contact structure such that all the submanifolds
$X_t\hook Y$ are isotropic. The contact line bundle $L$ on $Y$
is just the quotient $L=\tilde{\cal F}\times  {\Bbb C}/{\Bbb C}^*$
relative to the
natural multiplication map $\tilde{\cal F}\times{\Bbb C} \lon
\tilde{\cal F}\times  {\Bbb C}$, $(p,c)\rar (\lambda p, \lambda c)$, where
$\lambda\in {\Bbb C}^*$.
This can be summarized as follows.

\begin{proposition}\label{Iso}
Given an irreducible $G$-structure ${\cal G}\rar M$  with reductive
$G$. There is  canonically associated a complex contact manifold $(Y,L)$
containing
a $\dim M$-parameter family $\{X_t\hook Y\mid t\in M\}$  of isotropically
emdedded generalized flag varieties $X_t= G_s/P$, where
$G_s$ is the semisimple part of $G$ and $P$ is the parabolic subgroup
of $G_s$ leaving invariant a highest weight vector in the typical fibre of
$\Omega^1 M\rar M$  up to a scale.
\end{proposition}

Let $e$ be any point of $\tilde{\cal F}\subset \OM\setminus 0_{\OM} $.
Restricting a "lowering of indices" map $T_e(\OM)
\stackrel{\lrcorner\,\omega}
{\lon}\Omega^1_e(\OM)$ to the subspace ${\cal D}_e$, one obtains
an injective map
$$
0\lon {\cal D}_e\stackrel{\lrcorner\,\omega}{\lon} {\cal N}_e^*,
$$
where ${\cal N}_e^*$ is the fibre of the conormal bundle of
$\tilde{\cal F}\hook\OM\setminus 0_{\OM}$. Therefore, the rank of the
distribution ${\cal D}$ is equal at most to $\rank\, {\cal N}^* =\dim M
- \dim X_t - 1$.
It is easy to check that $\rank\, \cal D$ is maximal possible if
and only if ${\cal G}$ is involutive.
In this case the contact manifold $Y$ associated to ${\cal G}$ has
dimension
\begin{eqnarray*}
\dim Y & = & \dim \tilde{Y} - 1 = \dim \tilde{\cal F} - \rank {\cal D} - 1\\
   &=& (\dim M + \dim X_t + 1) - (\dim M - \dim X_t - 1) - 1
 = 2\dim X_t + 1
\end{eqnarray*}
which means that the associated complete family
$\{X_t\hook Y\mid t\in M\}$
is  an analytic family of compact Legendre submanifolds.
This argument partly explains the following result

\begin{theorem}
Let $G\subset GL(m,{\Bbb C})$ be an arbitrary connected
semisimple Lie subgroup whose decomposition into a locally direct
product of
simple groups does not contain any of the groups $H$ considered in
Theorem~\ref{Main2}(ii).
If ${\cal G}$ is any involutive $G\times {\Bbb C}^*$-structure on an
$m$-dimensional manifold $M$, then there exists a complex contact
manifold
$(Y,L)$ and a Legendre submanifold $X\hook Y$ with $X=G/P$ for some
parabolic
subgroup $P\subset G$ and with $L_X$ being very ample, such that,
at least
locally, $M$ is canonically isomorphic to the associated Legendre moduli
space and ${\cal G} = {\cal G}_{ind}$.
\end{theorem}

In the case of involutive $G$-structures ${\cal G} \rar M$ with $G$ as in
Theorem~\ref{Main2}(i) one can still canonically identify the
base manifold $M$ with a Legendre moduli space, but now
${\cal G}$ is properly contained in ${\cal G}_{ind}$.

In conclusion, the earlier posed question ---
how large is the family of $G$-structures which can be constructed by
twistor methods of Theorem~\ref{Main}? --- has the following answer:
this family consists of involutive $G$-structures.

\paragraph{8. Affine connections with "very little torsion".}
With any  irreducible $G$-structure  ${\cal G}$ on a complex manifold $M$
one can associate the {\em torsion number}\, defined as follows
$$
l = \frac{1}{2} \left(\dim M - \dim G/P - \rank {\cal D}  -1 \right).
$$
Here $P$ is the parabolic subgroup of $G$ leaving invariant
up to a scale a highest weight vector in the typical fibre of $TM$,
and ${\cal D}$ is the distribution associated
to  ${\cal G}$ as explained in section 7.
We see that the torsion number $l$ is composed of two
very different parts:
the first part, $\dim M - \dim G/P $, encodes only "linear" information
about the particular irreducible representation of $G$,
while the second part, $\rank {\cal D} - 1$, measures how this  particular
representation is "attached" to the base manifold $M$.
It is not difficult to prove  that $l$ is always a
 non-negative {\em integer}.
This fact alone shows that the  proposed combination $l$
of four natural numbers
does give some insight into the structure of ${\cal G}$.
This impression can be further strengthened by the fact that $l$
has a nice geometric interpretation. Remember that, by
Proposition~\ref{Iso}, the $G$-structure $\cal G$ gives rise to
a complex contact manifold $(Y,L)$ and a family $\{X_t\hook Y\mid t\in M\}$
of isotropic submanifolds parameterised by $M$. Then, in these terms,
$$
l=  \frac{1}{2}(\dim Y - 1) - \dim X_t ,
$$
 i.e.\ $l$ measures how much
$X_t$ lacks to be a Legendre submanifold.

Why is $l$ called a {\em torsion}\, number?  It is not difficult to
show that the torsion number of any {\em 1-flat }\, $G$-structure is zero.
Therefore, a $G$-structure on $M$ may have non-vanishing $l$  only if
it has a non-vanishing invariant torsion (but not vise versa, as we shall
see in a moment). Moreover, the larger $l$ is, the less integrable is
the distribution ${\cal D}$ and, in this sense, the "larger" is the
invariant torsion.

\begin{definition}
The torsion number of an affine connection  $\nabla$ is the torsion
number of the associated holonomy bundle ${\cal G}_{\nabla}$.
\end{definition}

\begin{definition}
An affine connection with torsion is said to have  very little torsion
if its torsion number is zero.
\end{definition}

The class of affine connections  with very little torsion  is a sibling
of the class torsion-free affine connections  in the sense that both
these classes (and only these two classes) have involutive holonomy
bundles (it is not difficult to show that a $G$-structure has
$l=0$ if and only if it is involutive). This means  in particular
that {\em all}\, connections
of both types can be constructed  by twistor methods on appropriate
Legendre moduli spaces.  Another conclusion is that
the class of affine connections with very little torsion is non empty
--- one can construct plenty of them using Legendre deformations
problems $"X\hook (Y,L)$" such that $H^1(X,L_X\otimes S^2(J^1 L_X)^*)
\neq 0$. This does not mean, however, that this class is enormously
large --- on the contrary, as the table below shows,
the  list of irreducibly acting holonomies of
affine connections with very little torsion must be very restricted.
\vspace{8 mm}

\begin{tabular}{|c||c|c|}
\hline
\multicolumn{1}{|c||}{ }  &  \multicolumn{2}{c|}{ The only possible
            holonomy representations  of $G$} \\
\multicolumn{1}{|c||}{Group }  &  \multicolumn{2}{c|}{in the class of
            connections with ``very little torsion''}\\
\cline{2-3}
$G$ & \ \ \ \ \ \ \ \ representation \ \ \ \ \ \ \ \  & $\dim$ \\
\hline
\hline
$SL(2,\Bbb C)$ & \begin{picture}(20,20)
                \put(6,2){\circle*{4}}
                \put(3,7){\shortstack{$k$}}
                \end{picture}
                 $k\geq 4$ & $k(k+1)/2$ \\
\hline

$SL(2,{\Bbb C})SL(2,\Bbb C)$ &
                              \begin{picture}(1,20)
                              \put(-2,2){\circle*{4}}
                              \put(-4,7){\shortstack{$1$}}
                              \end{picture}
                              $\otimes$
                             \begin{picture}(9,20)
                             \put(1,2){\circle*{4}}
                             \put(-1,7){\shortstack{$k$}} \,
                             \end{picture}\ \ $k\geq 2$ & $2k+2$  \\
                           &
                             \begin{picture}(1,20)
                             \put(-2,2){\circle*{4}}
                             \put(-4,7){\shortstack{$1$}}
                             \end{picture}
                             $\otimes$
                             \begin{picture}(9,20)
                             \put(1,2){\circle*{4}}
                             \put(-1,7){\shortstack{$k$}} \,
                             \end{picture}\ \  $k\geq 2$ & $3k+3$  \\
\hline

$SL(3,\Bbb C)$ & \begin{picture}(30,20)
                 \put(0,2){\line(5,0){30}}
                 \put(0,2){\circle*{4}}
                 \put(30,2){\circle*{4}}
                 \put(-2,6){\shortstack{1 }}
                 \put(28,6){\shortstack{1 }}
                 \end{picture}& 8  \\

               & \begin{picture}(30,20)
                 \put(0,2){\line(5,0){30}}
                 \put(0,2){\circle*{4}}
                 \put(30,2){\circle*{4}}
                 \put(-2,6){\shortstack{1 }}
                 \put(28,6){\shortstack{2 }}
                 \end{picture}& 15 \\
\hline

$Sp(4,\Bbb C)$ & \begin{picture}(30,20)
                 \put(0,0.7){\line(5,0){30}}
                 \put(0,3.3){\line(5,0){30}}
                 \put(0,2){\circle*{4}}
                 \put(30,2){\circle*{4}}
                 \put(12,6){\line(5,-2){10}}
                 \put(12,-2.5){\line(5,2){10}}
                 \put(-2,6){\shortstack{1 }}
                 \put(30,6){\shortstack{1 }}
                 \end{picture} & 16 \\

               & \begin{picture}(30,20)
                 \put(0,0.7){\line(5,0){30}}
                 \put(0,3.3){\line(5,0){30}}
                 \put(0,2){\circle*{4}}
                 \put(30,2){\circle*{4}}
                 \put(12,6){\line(5,-2){10}}
                 \put(12,-2.5){\line(5,2){10}}
                 \put(-2,6){\shortstack{2 }}
                 \put(30,6){\shortstack{0 }}
                 \end{picture} & 14 \\

               & \begin{picture}(30,20)
                 \put(0,0.7){\line(5,0){30}}
                 \put(0,3.3){\line(5,0){30}}
                 \put(0,2){\circle*{4}}
                 \put(30,2){\circle*{4}}
                 \put(12,6){\line(5,-2){10}}
                 \put(12,-2.5){\line(5,2){10}}
                 \put(-2,6){\shortstack{3 }}
                 \put(30,6){\shortstack{0 }}
                 \end{picture} & 30 \\
\hline

$G_2$  &  \begin{picture}(30,20)
          \put(0,0){\line(5,0){30}}
          \put(0,2){\line(5,0){30}}
          \put(0,4){\line(5,0){30}}
          \put(0,2){\circle*{4}}
          \put(30,2){\circle*{4}}
          \put(12,6){\line(5,-2){10}}
          \put(12,-2.5){\line(5,2){10}}
          \put(-2,6){\shortstack{0 }}
          \put(28,6){\shortstack{2 }}
          \end{picture}
          & 27 \\

       &  \begin{picture}(30,20)
          \put(0,0){\line(5,0){30}}
          \put(0,2){\line(5,0){30}}
          \put(0,4){\line(5,0){30}}
          \put(0,2){\circle*{4}}
          \put(30,2){\circle*{4}}
          \put(12,6){\line(5,-2){10}}
          \put(12,-2.5){\line(5,2){10}}
          \put(-2,6){\shortstack{0 }}
          \put(28,6){\shortstack{3 }}
          \end{picture}
          & 77 \\
\hline
\end{tabular}
\vspace{7 mm}

\noindent  Here irreducible representations
are written in the notations of \cite{BE}).

In conclusion we note that the problem of classifying all irreducibly
acting reductive holonomies of  affine connections with zero or very
little torsion has a strong purely symplectic flavour --- it
is nearly equivalent to the problem of classifying all
generalized flag varieties $X$ which can be realized as  complex
 Legendre submanifolds of  {\em non-trivial }\, contact manifolds $Y$
(non-trivial in the sense that the germ of $Y$ at $X$ is {\em not}\,
isomorphic to the germ of the total space of the jet bundle $J^1L_X$,
for some line bundle $L_X\rar X$, at its zero section).

\end{document}